\documentclass [12pt] {article}
\usepackage{amsfonts}
\usepackage{amssymb}

\newcommand{\qed} {\hspace {0.1in} \rule {1.5mm} {3.5mm}}

\newtheorem{lemma}{Lemma}[section]

\newtheorem{theorem}{Theorem}

\newtheorem{proposition}{Proposition}[section]

\def\lnt{\lim_{n\to\infty}}

\def\limn{\lim_{n\to\infty}}
\def\gseq{\{G_n\}_{n=1}^\infty}

\def\<{\langle}
\def\>{\rangle}

\def\proof{\smallskip\noindent{\it Proof.} }

\def\deg{\mbox{deg}\,}

\def\cA{\mbox{$\cal A$}}
\def\Mar{M(\cA_r)}
\def\Mbr{M(\cB_{r+1})}
\def\cB{\mbox{$\cal B$}}
\def\cC{\mbox{$\cal C$}}

\def\cG{\mbox{$\cal G$}}
\def\cU{\mbox{$\cal U$}}
\def\cV{\mbox{$\cal V$}}
\def\cW{\mbox{$\cal W$}}
\def\urd{\cU^{r,d}}
\def\vrd{\cV^{r,d}}
\def\wrdg{\cW^{r,d}_G}

\def\to{\rightarrow}

\begin{document}
\title{On limits of finite graphs} 
\author{\sc G\'abor Elek
\footnote {The Alfred Renyi Mathematical Institute of
the Hungarian Academy of Sciences, P.O. Box 127, H-1364 Budapest, Hungary.
email:elek@renyi.hu}}
\date{}
\maketitle
\vskip 0.2in
\noindent{\bf Abstract.}  We prove that for any weakly convergent
sequence of finite graphs with bounded vertex degrees,
there exists a topological limit graphing. 
\vskip 0.3in
\section{Introduction}

First let us recall the notion of  weak convergence of finite
graphs \cite{BS}. We need some definitions and notations. A {\it rooted graph}
is a simple graph with a distinguished vertex (the root). Two rooted graphs
are called rooted isomorphic if there exists a graph isomorphism between
them mapping one root to the other one. A rooted $r$-ball is a rooted graph
$G(V,E)$ such that $\sup_{y\in V(G)} d_G(x,y)\leq r$, where $x$ is the root of
$G$ and $d_G$ is the shortest path distance. For $d\geq 1$, we
denote by $\urd$ the finite
set of rooted isomorphism classes of rooted $r$-balls with vertex degrees
bounded by $d$. If $G(V,E)$ is a graph with $\sup_{y\in V(G)} \deg(y)\leq d$
and $\cA\in \urd$, then $T(G,\cA)$ denotes the set of vertices $v$ in $V(G)$
such that $\cA$ represents the rooted isomorphism class of the
$r$-neighborhood of $v$, $B_r(v)$. Set $p_G(\cA):=\frac{|T(G,\cA)|}{|V(G)|}\,.$

That is $G$ determines a probability distribution on $\urd$ for any $r\geq 1$.
Now let $\{G_n\}_{n=1}^\infty$ be a sequence of finite simple connected
graphs with
vertex degrees bounded by $d$ and $\lnt |V(G_n)|=\infty$. We say that $\gseq$
is {\it weakly convergent} if for any $r\geq 1$ and $\cA\in\urd$,
$\lnt p_{G_n}(\cA)$ exists.

\noindent
Now let $X$ be a standard Borel space with a measure $\mu$. Suppose that
$T_1, T_2,\dots, T_k$ are measure preserving Borel involutions of $X$. Then
the system $\cG=\{X,T_1,T_2,\dots,T_k,\mu\}$ is called a {\it measurable
graphing}. A measurable graphing $\cG$ determines an equivalence relation on
the points of $X$. Simply, $x\sim_{\cG} y$ if there exists a sequence of
points
$\{x_1, x_2,\dots, x_m\}\subset X$ such that
\begin{itemize}
\item $x_1=x, x_m=y$
\item $x_{i+1}=T_j(x_i)$ for some $1\leq j \leq k$.
\end{itemize}

Thus there exist natural simple graph structures on the equivalence classes,
the leafgraphs. Here $x$ is adjacent to $y$, if $x\neq y$ and $T_j(x)=y$
for some $1\leq j \leq k$. Now if $\cA\in\urd$, we denote by $p_{\cG}(\cA)$
the $\mu$-measure of the points $x$ in $X$ such that the rooted
$r$-neighborhood of $x$ in $\cA$ in its leafgraph is isomorphic to $\cA$.
We say that $\cG$ is the limit graphing of the weakly convergent sequence
$\gseq$, if for any $r\geq 1$ and $\cA\in\urd$
$$ \limn p_{G_n}(\cA)= p_{\cG} (\cA)\,.$$
Combining the ideas of \cite{BS} and \cite{Gab} one can see that
for any weakly convergent graph sequence there exists a measurable limit
graphing $\cG$. Note that for dense graph sequences an analogous construction\
is given in \cite{Lovasz}.

\noindent
If $X$ is a compact metric space with a Borel measure $\mu$ and
$T_1,T_2,\dots, T_k$ are continuous measure preserving involutions
of $X$, then $\cG=\{X,T_1,T_2,\dots, T_k,\mu\}$ is a {\it topological
graphing}. The goal of this note is to give a simple, self-contained proof
of the following theorem.
\begin{theorem}
\label{t1}
If $\gseq$ is a weakly convergent system of finite connected graphs, then
there exists a topological graphing $\cG$ such that for any $r\geq 1$ and
$\cA\in\urd$, $\limn p_{G_n}(\cA)= p_{\cG} (\cA)\,.$
\end{theorem}

\section{The construction of the limit graphing}
Our graphing construction is motivated by the ideas
of  \cite{BS},\cite{AL} and \cite{Gab}.
Let $\gseq$ be a weakly convergent sequence of finite connected graphs,
such that $|V(G_n)|\to \infty$ and $\sup_{v\in V(G_n)} \deg(y)\leq d$
for all $n\geq 1$.
Let $S$ be a set of $d+1$ elements. Color the edges of each $G_n$ by the
elements of $S$ to form an edge-coloring, that is if two edges have
a joint vertex then they have different colors. This can be done by
Vizing's Theorem. Now pick a sequence of integers $r_1\leq r_2 \leq \dots$
with the following property: Any graph $G$ with vertex degree bound $d$
has a vertex coloring by $r_i$ colors such that if $d_G(x,y)\leq i$, then
$x$ and $y$ are colored differently. Obviously if $r_i>(d+1)^i$, then these
conditions are satisfied. Now for any $i\geq 1$ and $n\geq 1$ we fix
a vertex coloring of $G_n$ by the elements of a set $Q_i$, $|Q_i|=r_i$,
satisfying the coloring conditions above, that is if $d_G(x,y)\leq i$, then
$x$ and $y$ are colored differently. We need a refined version of
rooted graph isomorphism. Let $G$ and $H$ be two rooted $r$-balls
with vertex degree bound $d$ equipped with an edge-coloring by $S$ and
a vertex-coloring with $\prod^r_{i=1} Q_i$ (satisfying the appropriate
 coloring conditions
for each $i\geq 1$.) Then $G$ and $H$ are called colored-isomorphic
if there exists a graph isomorphism $\phi:G\to H$, mapping root to root
and preserving both the edge and the vertex colorings. Let $\vrd$ be the
finite set of colored isomorphism classes.

\noindent
Again, for each $G_n$ and $r\geq 1$ one obtains a probability distribution
on $\vrd$. Let $ p^c_{G_n}(\cA)$ be the
probability of a colored type $\cA\in\vrd$.
 We may suppose that for any $\cA\in\vrd$, $\limn p^c_{G_n}(\cA)$ exists
(otherwise we can pick an appropriate subsequence). Let $\wrdg\subseteq\vrd$
be the set of colored types $\cA$ such that $\limn p^c_{G_n}(\cA)\neq 0$.
Let $\cA_{r-1}\in \cW^{r-1,d}_G, \cA_{r}\in \cW^{r,d}_G$, then
let $\cA_{r-1}\prec \cA_{r}$ if the $r-1$-ball around the root of
$\cA_{r}$ is just $\cA_{r-1}$. Since the vertices of $\cA_{r}$ are
colored by $\prod^r_{i=1}Q_i$ we can consider the restricted
$\prod^{r-1}_{i=1}Q_i$-coloring on the $r-1$-ball around the root, thus the
definition of $\prec$ is meaningful. 
We call an infinite sequence of classes $\cA_1\prec \cA_2 \prec \cA_3 \prec\dots$ a chain.
As usual the distance of two chains $\{\cA_i\}^\infty_{i=1}$
and $\{\cB_i\}^\infty_{i=1}$ is defined as $2^{-r}$, where $r$ is the minimal
index for which $\cA_r\neq\cB_r$. Denote the space of chains by $X$.
$X$ is clearly a compact metric space. If $a\in S$, then $T_a:X\to X$ is
defined
the following way. Let $x=\{\cA_1\prec \cA_2 \prec\dots\}\in X$. If there
is no edge colored by $a$, incident to the root of $x$, then let
$T_a(x)=x$. If there exists such an edge, then the root of $\cA_r$ is
connected with this $a$-colored edge
to a unique vertex $p$ of $\cA_r$ with $r-1$-ball $\cB_{r-1}$.
Clearly we have a chain $y=\{\cB_1\prec \cB_2 \prec\dots\}\,\in X.$ Define
$T_a(x)=y$. Then $T_a$ is a continuous involution.

\noindent
Now we define the measure $\mu$. Let $\Mar$ be the open-closed set of chains
starting with a
sequence $\cA_1\prec \cA_2 \prec\dots\prec\cA_r,$\, $\cA_r\in \wrdg$. Then let
$\mu(\Mar)=\limn  p^c_{G_n}(\cA_r)$. This way we define a Borel probability
measure on $X$.
\begin{proposition}
\label{p13}
$\mu$ is an invariant measure.
\end{proposition}
\proof
It is enough to prove, that if $M$ is an open-closed set as above and $a\in
S$, then $\mu(M)=\mu(T_a(M))$.
Again, if $\cA_r\in \wrdg$, then $\tau(G_n,\cA_r)$ denotes the number of
vertices $x$ of $G_n$ such that the $r$-ball around $x$ has colored type
$\cA_r$. Note that 
\begin{equation} \label {e4a}
T_a(\Mar)=\bigcup^{.}_{\cB_{r+1}\mid\,\cA_r\sim a} \Mbr\,,
\end{equation}
where
$\cB_{r+1}\mid\,\cA_r\sim a$ means that the summation
is taken over such $\cB_{r+1}$ that
\begin{itemize}
\item either there is no edge colored by $a$ going out from the root 
$x$ of $\cB_{r+1}$
and the $r$-ball around $x$ is of type $\cA_r$
\item or there exists $(x,y)\in E(G_n)$ colored by $a$ and
the $r$-ball around $y$ is of type $\cA_r$.
\end{itemize}
Also note that
\begin{equation} \label {e4b}
\tau(G_n,\cA_r)=\sum_{\cB_{r+1}\mid\,\cA_r\sim a}
\tau(G_n,\cB_{r+1})\,.
\end{equation}
Then 
$$\mu(T_a(\Mar))=\sum_{\cB_{r+1}\mid\,\cA_r\sim a} \mu(\Mbr)= $$
$$=\sum_{\cB_{r+1}\mid\,\cA_r\sim a}\limn p^c_{G_n}(\cB_{r+1})=
\limn \sum_{\cB_{r+1}\mid\,\cA_r\sim a} p^c_{G_n}(\cB_{r+1})= $$
$$=\limn \sum_{\cB_{r+1}\mid\,\cA_r\sim a} \frac {\tau(G_n, \cB_{r+1})}
{|V(G_n)|}=\limn \frac {\tau(G_n, \cA_{r})}
{|V(G_n)|}=\limn  p^c_{G_n}(\cA_{r})=\mu(\Mar)\,.
$$
Hence the proposition follows. \qed

\noindent
In order to finish the proof of Theorem \ref{t1} we only need to prove the 
following proposition.
\begin{proposition}
\label{p15}
If $x=\{\cA_1\prec \cA_2\prec\dots\}\in X$ then the $r$-ball around
$x$ in the leafgraph of $x$ has the same rooted isomorphism
class as $\cA_r$.
\end{proposition}
The proof shall be given in two lemmas.
\begin{lemma}
\label{l15}
Let $x=\{\cA_1\prec \cA_2\prec\dots\}\in X$ be a chain. Let
$\{i_1, i_2,\dots,i_s\}$ and \\ $\{j_1, j_2,\dots,j_t\}$ sequences
of elements of $S$ such that $s,t\leq r$. Suppose that 
$T_{i_1}T_{i_2}\dots T_{i_s}(v)= T_{j_1}T_{j_2}\dots T_{j_t}(v)$, where
$v$ is the root of $\cA_r$. Then
$T_{i_1}T_{i_2}\dots T_{i_s}(x)= T_{j_1}T_{j_2}\dots T_{j_t}(x)$. Note
that the involutions are defined both on the graphs and on the space $X$.
\end{lemma}
\proof
For any $r'>r$, if $\cA_r\prec \cA_{r'}$, then 
$T_{i_1}T_{i_2}\dots T_{i_s}(v)= T_{j_1}T_{j_2}\dots T_{j_t}(v)$
holds for $\cA_{r'}$. Thus
 $T_{i_1}T_{i_2}\dots T_{i_s}(x)= T_{j_1}T_{j_2}\dots T_{j_t}(x)$. \qed

\begin{lemma}
\label{l16}
Let $x$ be as above and $T_{i_1}T_{i_2}\dots T_{i_s}(v)\neq
T_{j_1}T_{j_2}\dots T_{j_t}(v)$.\\ Then 
$T_{i_1}T_{i_2}\dots T_{i_s}(x)\neq T_{j_1}T_{j_2}\dots T_{j_t}(x)$.
\end{lemma}
\proof
This is the point where we need the vertex coloring.
Recall that in $\cA_{2r}$, \\ $T_{i_1}T_{i_2}\dots T_{i_s}(v)$ and
$T_{j_1}T_{j_2}\dots T_{j_t}(v)$ are colored by $\prod^{2r}_{i=1} Q_i$ and
their $Q_{2r}$-colors are different. Hence if 
$$T_{i_1}T_{i_2}\dots T_{i_s}(x)=\{\cB_1\prec\cB_2\prec\dots\prec\cB_{2r}\prec
\dots\}$$
and
$$T_{j_1}T_{j_2}\dots T_{j_t}(x)=\{\cC_1\prec\cC_2\prec\dots\prec\cC_{2r}\prec
\dots\}$$
then $\cB_{2r}\neq \cC_{2r}$. Therefore
$T_{i_1}T_{i_2}\dots T_{i_s}(x)\neq T_{j_1}T_{j_2}\dots T_{j_t}(x)$. \qed

\noindent
Again let $x=\{\cA_1\prec \cA_2\prec\dots\}$ and $v$ be the root of $\cA_r$.
Suppose that $w\in \cA_r$ and $w=T_{i_1}T_{i_2}\dots T_{i_s}(v)$. Then
let $\phi(w):=T_{i_1}T_{i_2}\dots T_{i_s}(x)$. By Lemma \ref{l15},
$\phi$ is well-defined and by Lemma \ref{l16}, $\phi$ is injective.
Since any element of the $r$-ball around $x$ in the leaf-graph of $x$ is
in the form $T_{j_1}T_{j_2}\dots T_{j_t}(x)$, for some $t\leq r$ and
$\{j_1, j_2,\dots,j_t\}$, $\phi$ is also surjective. Hence $\phi$
is a bijection between $\cA_r$ and the $r$-ball around $x$ in
the leaf graph of $x$. Note that $\phi$ preserves adjacency. Thus
our proposition and consequently Theorem \ref{t1} follow.
 \qed

\end{document}